\documentclass[11pt,a4paper]{article}
\usepackage[francais,english]{babel}
\usepackage[utf-8]{inputenc}
\usepackage{hyperref}
\usepackage{latexsym,amssymb,amsmath,amsfonts,delarray,epsfig,eucal,graphicx,hyperref}
\usepackage{amsmath}
\usepackage{amssymb}
\usepackage[T1]{fontenc}
\usepackage{multicol}
\usepackage{calc}
\usepackage[all]{xy}
\usepackage{amsthm}

\theoremstyle{plain}
\newtheorem{theoreme}{\sc{Theorem}}[section]
\newtheorem{proposition}[theoreme]{Proposition}
\newtheorem{corollaire}[theoreme]{Corollary}
\newtheorem{lemme}[theoreme]{Lemma}

\theoremstyle{definition}
\newtheorem{definition}[theoreme]{Definition}

\theoremstyle{remark}
\newtheorem{remarques}[theoreme]{Remarks}

\newtheorem{remarque}[theoreme]{Remark}

\title{Smoothness in Relative Geometry}
\author{Florian Marty\\fmarty@math.ups-tlse.fr}
\date{Université Toulouse III - Laboratoire Emile Picard}
\begin{document}
\small
\maketitle
\setcounter{secnumdepth}{-1}
\section{Abstract}
\footnotesize{In \cite{tva}, Bertrand Toën and Michel Vaquié defined a scheme theory for a closed monoidal category $(\mathcal{C},\otimes,1)$. In this article, we define a notion of smoothness in this relative (and not necessarily additive) context which generalize the notion of smoothness in the category of rings. This generalisation consists practically in changing homological finiteness conditions into homotopical ones using Dold-Kahn correspondence. To do this, we provide the category $s\mathcal{C}$ of simplicial objects in a monoidal category and all the categories $sA-mod$, $sA-alg$ ($A\in scomm(\mathcal{C})$) with compatible model structures using the work of Rezk in \cite{r}. We give then a general notions of smoothness in $sComm(\mathcal{C})$. We prove that this notion is a generalisation of the notion of smooth morphism in the category of rings, is stable under compositions and homotopic pushouts and provide some examples of smooth morphisms in $\mathbb{N}-alg$ and $Comm(Set)$.} 
\setcounter{secnumdepth}{2}
\small
\tableofcontents
\setcounter{secnumdepth}{-1}
\section{Introduction}
In \cite{tva}, Bertrand Toën and Michel Vaquié defined a scheme theory for a closed monoidal category $(\mathcal{C},\otimes,1)$. In this article, we define a notion of smoothness in this relative context which generalize the notion of smoothness in the category of rings. The motivations for this work are that interesting objects in the non additive contexts $\mathcal{C}=Ens$ or $\mathbb{N}-mod$ are expected not to be schemes but Stacks. A theorem asserts for $\mathcal{C}=\mathbb{Z}-mod$ that quotients of schemes by smooth group schemes are in fact algebraic stacks. A near theorem is expected in a relative context. The first step is to get a definition for smooth morphism.\\
The following theorem give the good definition of smoothness that can be generalised to the relative context:
\begin{theoreme}\label{l}
Assume $\mathcal{C}=\mathbb{Z}-mod$ A morphism of rings $A\rightarrow B$ is smooth if and only if
\begin{enumerate}
\item The ring  $B$ is finitely presented in $A-alg$.
\item The morphism $A\rightarrow B$ is flat.
\item The ring $B$ is a perfect complex of $B\otimes_{A}B$ modules.
\end{enumerate}
\end{theoreme}
The flatness of $A\rightarrow B$ is in fact equivalent to $Tordim_{A}B=0$ hence the two last conditions are homological finiteness conditions. By the correspondence of Dold-Kan, the second condition can then be traduced in an homotopical condition. Finally, a result from \cite{tv} asserts that $B$ is a perfect complex in $ch(B\otimes_{A}B)$ if and only if it is homotopically finitely presentable in $sB\otimes_{A}B-mod$.\\ 
We provide then the category $s\mathcal{C}$ of simplicial objects in a monoidal category and all the categories $sA-mod$, $sA-alg$ with  model structures using the work of Rezk in \cite{r}. The classical functors between the categories $A-mod,A-alg,\;A\in sComm(\mathcal{C})$ induce Quillen functors between the corresponding simplicial categories. We give the following general definition for smoothness 
\begin{definition}
Let $A$ be in $sComm(\mathcal{C})$, a morphism $B\rightarrow C$ in $sA-alg$ is smooth if and only if
\begin{enumerate}
\item The simplical algebra $C$ is homotopically finitely presented in $sB-alg$.
\item The simplicial algebra $C$ has Tor dimension $0$ on $B$.
\item The morphism $C\otimes^{h}_{B}C\rightarrow C$ is homotopically finitely presented in $sC\otimes^{h}_{B}C-mod$.
\end{enumerate}
\end{definition}
The first condition imply the first condition of \ref{l} (\cite{tv}, 2.2.2.4) and there is equivalence for smooth morphisms of rings. When  $A,B,C$ are just rings, the Tor dimension $0$ imply that the derived tensor product is weakly equivalent to the tensor product. The equivalence with previous theorem for rings is then clear.\vskip 2pt
We prove that relative smooth morphisms are stable under composition and pushouts of Algebras. We finally provide examples of smooth morphism in relative non-additive contexts , for $\mathcal{C}=\mathbb{N}-mod$ or $\mathcal{C}=Ens$. In particular the affine lines $\mathbb{F}_{1}\rightarrow \mathbb{N}$ and $\mathbb{N}\rightarrow \mathbb{N}[X]$ are smooth respectively in $sComm(Set)$ and $s\mathbb{N}-mod$.

\section{Preliminaries}
\noindent{Let $(\mathcal{C},\otimes,1)$ be a complete and cocomplete closed symmetric monoidal category. In the category $\mathcal{C}$, there exists a notion of commutative monoid and for a given commutative monoid $A$, of $A$-module. Let $Comm(\mathcal{C})$ denotes the category of commutative monoids (with unity) in $\mathcal{C}$ and for $A\in Comm(\mathcal{C})$, $A-mod$ denotes the category of $A$-modules. It is well known that the category $A-mod$ is a closed monoidal tensored and cotensored category, complete and cocomplete. The category $Comm(A-mod)$ will be denoted by $A-alg$ and is described by the equivalence $A/Comm(\mathcal{C})\backsimeq A-alg$. A pushout in $A-alg$ is a tensor product in the sense that for commutative monoids $B,C\in A-alg$ $B\otimes_{A}C\backsimeq B\coprod_{A}C$.\vskip 2pt
\textbf{All along this work}, $(\mathcal{C},\otimes,1)$ is a locally finitely presentable monoidal category i.e. verify that the (Yoneda) functor $i:\mathcal{C}\rightarrow Pr(\mathcal{C}_{0})$, where $\mathcal{C}_{0}$ is the full subcategory of finitely presented objects, is fully faithful, that $\mathcal{C}_{0}$ is stable under tensor product ans contains the unity. The functor $Hom_{\mathcal{C}}(1,-)$, denoted $(-)_{0}:X\rightarrow X_{0}$ is called ''underlying set functor''. For $k\in \mathcal{C}_{0}$, the functors $Hom_{\mathcal{C}}(k,-)$, denoted $(-)_{k}:X\rightarrow X_{k}$ are called ''weak underlying set functor''. It is known that if $\mathcal{C}$ is a locally finitely presentable monoidal category, so are its categories of modules $A-mod$, $A\in sComm(\mathcal{C})$\vskip 2pt
Let $J$ denotes the set of isomorphism classes of the objects of $\mathcal{C}_{0}$. We assume that the adjunction $\xymatrix{\mathcal{C} _{(resp\;Comm(\mathcal{C}))}\ar@<1ex>[r]&Ens^{J}\ar@<1ex>[l]}$  induced is monadic. Thus, on simplicial categories, the theorem of Rezk \ref{rezk} will provide a (good) model structure. Moreover, as the forgetful functor from $A-mod$ to $\mathcal{C}$ preserves colimits for a monoid $A$ the functor from $A-mod$ to $Ens^{J}$ is also monadic. This is a consequence of the characterisation of monadic functors (see ).
\vskip 2pt
There are two fundamental adjunctions:
\begin{center}
$\xymatrix{\mathcal{C}\ar@<1ex>[r]^-{(-\otimes A)}&A-mod\ar@<1ex>[l]^-{i}}$\;\;\;\;\;\;\;\;\;\;$\xymatrix{\mathcal{C}\ar@<1ex>[r]^-{L}&Comm(\mathcal{C})\ar@<1ex>[l]^-{i}}$
\end{center}
where the forgetful functor $i$ is a right adjoint and the functor ''free associated monoid'' $L$ is defined by $L(X):=\coprod_{n\in \mathbb{N}}X^{\otimes n}/S_{n}$. In these adjunctions, $\mathcal{C}$ can be replaced by $B-mod$ for $B\in Comm(\mathcal{C})$, and $L$ by $L_{B}$ defined by $L_{B}(M):=\coprod_{n\in \mathbb{N}}M^{\otimes_{B} n}/S_{n}$. Let $\varphi$ (resp $\varphi_{B}$) and $\psi$ (resp $\psi_{B}$) denote these adjunctions for the category $\mathcal{C}$ (resp $B-mod$). For $X\in \mathcal{C}$ and $M\in A-mod$, $\varphi:Hom_{\mathcal{C}}(X,M)\rightarrow Hom_{A-mod}(X\otimes A,M)$ is easy to describe :
\begin{center}
$\varphi: f\rightarrow \mu_{M}\circ Id_{A}\otimes f$\vskip 3pt$\varphi^{-1}:g\rightarrow g\circ (Id_{X}\otimes i_{A})\circ r_{X}^{-1}$
\end{center}
Let $s\mathcal{C}$ denotes the category of simplicial objects in $\mathcal{C}$. There is a functor ''constant simplicial object'' denoted $k$ from $\mathcal{C}$ to $s\mathcal{C}$ which is right adjoint to the functor $\pi_{0}$ from $s\mathcal{C}$ to $\mathcal{C}$ defined by $\pi_{0}(X):=Colim(\xymatrix{X[1]\ar@<1ex>[r]\ar@<-1ex>[r]&X[0]\ar[l]})$. The tensor product of $\mathcal{C}$ induces a tensor product on $s\mathcal{C}$, its unity element is $k(1)$. For $A$ in $Comm(\mathcal{C})$, $sA-mod$ and $sA-alg$ will denote respectively the simplicial categories $sk(A)-mod$ and $sk(A)-alg$. As $sComm(\mathcal{C})\backsimeq Comm(s\mathcal{C})$,  we will always refer to simplicial category of commutative monoids in $s\mathcal{C}$ as $sComm(\mathcal{C})$. The functor induced by $L$ on simplical categories will be denoted $sL$.
The functor $i:\mathcal{C}\rightarrow Pr(\mathcal{C}_{0})$ induces a functor $si:s\mathcal{C}\rightarrow sPr(s\mathcal{C}_{0})$. \vskip 2pt We need finally hypotheses to endow $s\mathcal{C}$, $sComm(\mathcal{C})$, and for $A\in sComm(\mathcal{C})$, $sA-mod$ and $sA-alg$ with compatible model structures. One solution of this question is to assume that the natural functors from $s\mathcal{C}$ and $sComm(\mathcal{C})$ to $sSet^{J}$ are monadic, where $J$ is the set of isomorphism classes of $J$. The characterisation of monadic functors of \cite{bc} implies that for any commutative simplicial monoid $A$, the induced functors from $sA-mod$ to $sSet^{J}$ is also monadic. 
\setcounter{secnumdepth}{2}
\section{General Theory}
\subsection{Simplicial Categories and Simplicial Theories}
\begin{definition}
A simplicial theory is a monad (on $sSet^{J}$) commuting with filtered colimits.
\end{definition}

\begin{theoreme}(Rezk)\label{rezk}\\
Let T be a simplicial theory in $sSet^{J}$, then $T-alg$ admits a simplicial model structure. $f$ is a Weak equivalence or a fibration in $T-alg$ if and only if so is its image in $sSet^{J}$ (for the projective model structure). Moreover, this Model category is right proper.  
\end{theoreme}

\begin{proposition}Model structures on the simplicial categories.
\begin{enumerate}
\item Let $A=(A_{p})$ be a commutative monoid in $s\mathcal{C}$. The monadic adjunctions $\xymatrix{A_{p}-mod\ar@<1ex>[r]&Set^{J}\ar@<1ex>[l]}$ induce a monadic adjunction $\xymatrix{sA-mod\ar@<1ex>[r]&sSet^{J}\ar@<1ex>[l]}$ i.e. there is an equivalence $sA-mod\backsimeq T_{A}-alg$ = where $T_{A}$ is the monad induced by adjunction. In particular $s\mathcal{C}\backsimeq T_{1}-alg$. 
\item Let $A=(A_{p})$ be a commutative monoid in $s\mathcal{C}$. The monadic adjunctions $\xymatrix{A_{p}-alg\ar@<1ex>[r]&Ens^{J}\ar@<1ex>[l]}$ induce a monadic adjunction $\xymatrix{sA-alg\ar@<1ex>[r]&\ar@<1ex>[l]sSet^{J}\ar@<1ex>[l]}$ i.e. there is an equivalence $sA-alg\backsimeq T^{c}_{A}-alg$ where $T^{c}_{A}$ is the monad induced by adjunction. In particular $sComm(\mathcal{C})\backsimeq T^{c}_{1}-alg$.
\end{enumerate}
\end{proposition}
Proof\\As explained in the preliminaries, this is due to the characterisation of monadic functors (\cite{bc}).
\vskip 2pt
$\hfill\blacklozenge$
\begin{remarque}
The right adjoints functors all commute with filtered colimits. So do the monads which are then simplicial theories on $sSet^{J}$.
\end{remarque}

\begin{corollaire}
Let $A$ be a commutative monoid in $s\mathcal{C}$. The categories $s\mathcal{C}$ and $sA-mod$ and $sComm(\mathcal{C})$ are Model categories. Moreover, the functors $(A\otimes-)$ and $sL$ are left Quillen and their adjoints preserve by construction weak equivalences and fibrations.
\end{corollaire}
\begin{theoreme}
The Category $s\mathcal{C}$ $(resp\;sA-mod)$ is a monoidal model category
\end{theoreme}
Proof:\\The proof for $s\mathcal{C}$ and $sA-mod$ are similar, so let us prove it for $s\mathcal{C}$. Let $I$, $I'$ be respectively the sets of generating cofibration and generating trivial cofibration. $I$ and $I'$ are the image by the left adjoint functor respectively of generating cofibration and generating trivial cofibration in $sSet^{J}$. We juste have to prove (cf \cite{h} chap IV) that $I\square I$ is a set of cofibrations and $I\square I'$ and $I'\square I$ are sets of trivial cofibrations. It is true for generating cofibration and generating trivial cofibration in $sSet^{J}$, which are all morphisms concentrated in a fixed level. Moreover, it is easy to verify that the functor $sK$ commutes with $\square$ of morphisms concentrated in one level. So it is true in $s\mathcal{C}$. The second axiom is clearly verified, as $1$ is cofibrant.
\vskip 2pt
$\hfill\blacklozenge$

\subsection{Compactly Generated Model categories}
\begin{definition}
Let $\mathcal{M}$ be a cofibrantly generated simplicial model category and $I$ be the set of generating cofibrations.
\begin{enumerate}
\item An object $X\in I-cell$ is strictly finite if and only if there exists a finite sequence 
\begin{center}
$\xymatrix{\emptyset=X_{0}\ar[r]&X_{1}\ar[r]&...\ar[r]&X_{n}=X}$
\end{center}
and $\forall\;i$ a pushout diagram:
\begin{center}
$\xymatrix{X_{i}\ar[r]\ar[d]&X_{i+1}\ar[d]\\A\ar[r]_{u_{i}}&B}$
\end{center}
with $u_{i}\in I$.
\item An object $X\in I-cell$ is finite if and only if it is weakly equivalent to a strictly finite object.
\item An object $X$ is homotopically finitely presented if and only if for any filtered diagram $Y_{i}$, the morphism :
\begin{center}
$Hocolim_{i}Map(X,Y_{i})\rightarrow Map(X,Hocolim_{i}Y_{i})$
\end{center}
is an isomorphism in $Ho(sSet)$
\item A model category $\mathcal{M}$ is compactly generated if it is cellular, cofibrantly generated and if the domains and codomains of generating cofibration and generating trivial cofibration are cofibrant, $\omega$-compact and $\omega$-small (Relative to $\mathcal{M}$).
\end{enumerate}
\end{definition}
\begin{proposition}Let $\mathcal{M}$ be a compactly generated model category.
\begin{enumerate}
\item For any filtered diagram $X_{i}$, the natural morphism $Hocolim_{i}X_{i}\rightarrow Colim_{i}X_{i}$ is an isomorphism in $Ho(\mathcal{M})$.
\item Assume that filtered colimits are exact in $\mathcal{M}$. Then homotopically finitely presented objects in $\mathcal{M}$ are exactly objects equivalent to weak retracts of strictly finite $I-cell$ objects.
\end{enumerate}
\end{proposition}
\begin{proposition}\label{pja}
\begin{enumerate}
\item The simplicial model category $sSet^{J}$ is compactly generated.
\item The categories of simplicial algebras over a simplicial theory are compactly generated.  
\end{enumerate}
\end{proposition}
\begin{lemme}\label{pjaa}
Let $A$ be in $sComm(\mathcal{C})$. Let $u^{j}$ be the family of images by the left adjoint functor in $sA-mod$ (resp $sA-alg$) of elements $*_{j}$ of $sSet^{J}$ defined by $*$ on level $j$ and $\emptyset$ on other levels. Any codomains of a generating cofibrations of $sA-mod$ (resp $sA-alg$) is weakly equivalent to an object $u^{j}$.  Any domain of a generating cofibration is weakly equivalent, for a given element $j$ in $J$, to an object obtained from the initial object (denoted $\emptyset$) and $u^{j}$ in a finite number of homotopic pushouts.
\end{lemme}
Proof:\\Generating cofibrations of $sA-mod$ are images of generating cofibrations of $sSet^{J}$ by the left adjoint functor. Generating cofibrations of $sSet$ are morphisms $\delta\Delta^{p}\rightarrow \Delta^{p}$. Their codomain is contractible, thus so are the codomains of generating cofibrations of $sSet^{J}$ for the projective model structure, and their image by by the left adjoint is weakly equivalent to the unity $1$. For the domains, consider the relation $\delta\Delta^{p+1}\backsimeq \Delta^{p+1}\coprod_{\delta\Delta^{p}}^{h}\Delta^{p+1}\backsimeq *\coprod_{\delta\Delta^{p}}^{h}*$ and $\delta\Delta^{0}=\emptyset$. Domains of generating cofibration in $sSet^{J}$ for the projective model structure are objects $(\delta\Delta^{p,j})_{p\in N,\;j\in V}$ defined in level $i\neq j$ by $\emptyset$ and in level $j$ by $\delta\Delta^{p}$ and verify the relation
\begin{center}
$(\delta\Delta^{p,j}) \backsimeq  *_{j}\coprod_{(\delta\Delta^{p-1,j})}^{h}*_{j}$
\end{center}
Clearly $\delta\Delta^{0,j} = \emptyset$ and $\delta\Delta^{1,j}=*_{j}$. Let $u^{p,j}$ denote the image of $\delta\Delta^{p,j}$. For all $j$, $u^{p,j}$ is obtained in a finite number of pushouts from $\emptyset$ and $u^{j}$.
\vskip 2pt
$\hfill\blacklozenge$
\begin{corollaire} of proposition \ref{pja} and lemma \ref{pjaa}.
\begin{enumerate}
\item The Simplicial Model categories $s\mathcal{C}$, $sA-mod$ ($A\in sComm(\mathcal{C})$), $sComm(\mathcal{C})$ and $sA-alg$ ($a\in sComm(\mathcal{C})$) are compactly generated.
\item Homotopically finitely presented objects of $sA-mod$ $(resp sA-alg)$ are exactly objects weakly equivalent to weak retracts of strictly finite $I-Cell$ objects.
\item The sub-category of $Ho(ssA-mod)$ $(resp\;Ho(sA-alg))$ $Ho(sA-mod)_{c}\;(resp\;Ho(sA-alg)_{c})$ Consisting of homotopically finitely presented objects is the smallest full sub-category of $Ho(sA-mod)$ (resp $Ho(sA-alg)$) containing the family $(u^{1,j})_{j\in V}$ $(resp\;(u_{a}^{1,j})_{j\in V})$, and stable under retracts and homotopic pushouts.
\end{enumerate}
\end{corollaire}
Proof of $iii$:\\Let $\mathcal{D}$ be the smallest full sub-category of $Ho(s\mathcal{C})$ containing $(u^{j})_{j\in V}$ $(resp\;(u_{A}^{j})_{j\in V})$, the initial object $\emptyset$ and stable under retracts and homotopic pushouts. Clearly, by $ii$, as $(\emptyset\rightarrow u^{j})_{j\in V}$ are  generating cofibrations of $s\mathcal{C}$, $Ho(s\mathcal{C})_{c}$ contains the family $(u^{j})_{j\in V}$, and is stable under retracts and homotopic pushouts. Thus $\mathcal{D}\subset Ho(s\mathcal{C})_{c}$. Reciprocally, let $X$ be an object of $Ho(s\mathcal{C})_{c}$, by $ii$, $X$ is isomorphic to a weak retract of a strictly finite $I-cell$ object. Therefore, there exists $n$ and $X_{0}...X_{n}$ such that:
\begin{center}
$\xymatrix{\emptyset=X_{0}\ar[r]&X_{1}\ar[r]&...\ar[r]&X_{n}=X}$
\end{center}
and $\forall j \in \{0,..,n-1\}$, $\exists K\rightarrow L$, a generating cofibration such that:
\begin{center}
$\xymatrix{X_{j}\ar[r]&X_{j+1}\\K\ar[r]\ar[u]&L\ar[u]}$
\end{center}
is a pushout diagram. Now, as domains and codomains of generating cofibrations are in $\mathcal{D}$, $X$ is in $\mathcal{D}$. 
\vskip 2pt
$\hfill\blacklozenge$

\subsection{Categories of Modules and Algebras}
\begin{proposition}Let $A$ be in $sComm(\mathcal{C})$ and $B$ be a simplicial monoid in $sA-alg$, cofibrant in $sA-mod$ .
\begin{enumerate}
\item The forgetful functor from $sB-mod$ to $sA-mod$ preserves cofibrations
\item The forgetful functor from $sA-alg$ to $sA-mod$ preserves cofibrations whose domain is cofibrant in $sA-mod$. In particular, it preserves cofibrant objects.
\end{enumerate}
\end{proposition}
Proof\\In each case, we just have to prove it for generating cofibrations and then generalize it to any cofibration by the small object argument.\vskip 2pt Proof of $i$: First, we choose a generating cofibration in $sB-mod$. As generating cofibration of $sB-mod$ are images of generating cofibrations of $sSet^{J}$, we set $L\rightarrow M$, a generating cofibration in $sSet^{J}$. Let $K_{A}$, $K_{A}$ denotes respectively the left adjoint functors (from $sSet^{J}$) for $sA-mod$ and $sB-mod$. The axiom of stability under $\square$ implies that the morphism $(\emptyset\rightarrow B)\square (K_{A}(L)\rightarrow K_{A}(M))$ is a cofibration in $sA-mod$. This morphism is in fact $K_{B}(L)=B\otimes_{A} K_{A}(L)\rightarrow B\otimes_{A} K_{A}(M)=K_{B}(M)$, hence generating cofibrations of $sB-mod$ are cofibrations in $sA-mod$.\vskip 2pt Proof of $ii$: As for $i$, let $N\rightarrow M$ be a generating cofibration in $sSet^{J}$. Let $L^{s}$ denotes the functor '' free associated commutative monoid'' of $sSet^{J}$. The functors $L$ (resp $sL_{A}$ in $sA-mod$) and $K_{A}$ are defined by colimits and so commute up to isomorphisms. That means that $K_{A}\circ L^{s}\backsimeq sL_{A}\circ K_{A}$. So the generating cofibration of $sA-alg$ corresponding to $N\rightarrow M$ is isomorphic to $K_{A}(L(N))\rightarrow K_{A}(L(M))$. To prove that it is a cofibration in $sA-mod$, we have then to prove that the morphism $L(N)\rightarrow L(M)$ is injective levelwise and this is clear as for any morphism $N^{\otimes n}/S_{n}\rightarrow M^{\otimes n}/S_{n}$ is injective. Thus any generating cofibration of $sA-alg$ is a cofibration in $sA-mod$. In fact it is a generating cofibration of $sA-mod$. To use the small object argument (of $sA-alg$), we need to verify that it preserves cofibrations in $sA-mod$. In fact, we need to check that an homotopic pushout in $sA-alg$ of a cofibration in $sA-mod$ is still a cofibration in $sA-mod$. We let the reader verify that it is a consequence of the axiom of stability by $\square$. Finally, the forgetful functor preserves cofibrations and, as $A$ is cofibrant in $sA-mod$, any cofibrant object of $sA-alg$ is also cofibrant in $sA-mod$.
\vskip 2pt
$\hfill\blacklozenge$
\begin{lemme}
Let $A\rightarrow B\;\in sComm(\mathcal{C})$ be a trivial cofibration between cofibrant objects. The categories of module are equivalent i.e. $Ho(sA-mod)\backsimeq Ho(sB-mod)$.
\end{lemme}
\noindent{Proof:\\We must prove that for $X$ cofibrant in $sA-mod$ and $y$ fibrant in $sB-mod$, $\varphi_{a}(f):x\otimes_{A}B\rightarrow y$ is a weak equivalence in $sB-mod$ if and only if so is $f: X\rightarrow Y$ in $sA-mod$. By previous lemma, $A\rightarrow B$ is a trivial cofibration in $sA-mod$. Thus as $X$ is cofibrant, using the axiom of stability under $\square$, $g:X\rightarrow B\otimes_{A} X$ is a weak equivalence in $sA-mod$. By construction of the adjunction $\varphi_{A}$, the following diagram is commutative :
\begin{center}
$\xymatrix{X\ar[r]\ar@/^-1pc/[rrr]_{f}&X\otimes_{A}B\ar[r]\ar@/^1pc/[rr]^{\varphi_{A}(f)}&Y\otimes_{A}B\ar[r]&Y}$
\end{center} 
Thus $f=g\circ\varphi_{A}(f)$. Finally, $\varphi_{A}(f)$ is a weak equivalence in $sA-mod$ if and only if so it is in $sB-mod$ and the \textit{two out of three} axiom ends the proof.
\vskip 2pt
$\hfill\blacklozenge$
\begin{proposition}
Let $f:A\rightarrow B\;\in sComm(\mathcal{C})$ be a weak equivalence between cofibrant objects. The categories of module are equivalent ie $Ho(sA-mod)\backsimeq Ho(sB-mod)$.
\end{proposition}
Let $r_{c}$ be the fibrant replacement of $sComm(\mathcal{C})$, then by previous lemma, the homotopical categories of modules over $A$ and $r_{c}A$ (resp $B$ and $r_{c}B$) are equivalent. Thus $A$ and $B$ can be taken fibrant and $f$ is an homotopy equivalence i.e. $\exists$ $g$ such that $f\circ g$ and $g\circ f$ are homotopic to identity. The following diagrams are commutative:
\begin{center}
$\xymatrix{B\ar[rd]^{Id}\ar[d]_{i_{0}}\\B_{1}\ar[r]^{h}&B\\B\ar[u]^{i_{1}}\ar[ru]_{f\circ g}}$\;\;\;\;\;\;\;\;\;\;\;\;\;\;\;\;\;\;
$\xymatrix{Ho(B-mod)\\Ho(B_{1}-mod)\ar[u]_{i_{0}^{*}}\ar[d]^{i_{1}^{*}}&Ho(B-mod)\ar[ld]_{(f\circ g)^{*}}\ar[l]^{h^{*}}\ar[lu]^{Id}\\Ho(B-mod)}$
\end{center}
where $i_{0}$ and $i_{1}$ are cofibrations and have the same right inverse  $p$ i.e. such that $p\circ i_{1}=p\circ i_{0}= Id_{B}$. the morphism $h$ is a trivial fibration thus $i_{0}$ is a weak equivalence. By previous lemma, $i_{0}^{*}$ is an equivalence of categories. Thus so is $p^{*}$. As $i_{1}^{*}$ and $i_{0}^{*}$ are both inverses of $p^{*}$, they are isomorphic and $i_{1}^{*}$ is also an equivalence. Finally, $h^{*}$ is an equivalence and so is $(f\circ g)^{*}$. The same method prove that $(g\circ f)^{*}$ is an equivalence. 
\vskip 2pt
$\hfill\blacklozenge$

\subsection{Finiteness Conditions}
\begin{definition}\label{fin}
Let $q_{c}$ be a cofibrant replacement in $sComm(\mathcal{C})$ and $f:A\rightarrow B$ be a morphism in $sComm(\mathcal{C})$.
\begin{list}{$\triangleright$}{}
\item The morphism $f$ is homotopically finite (denoted $hf$) if $B$ is homotopically finitely presented in $sq_{c}A-mod$. 
\item The morphism $f$ is homotopically finitely presented (denoted $hfp$) if $B$ is homotopically finitely presented in $sq_{c}A-alg$.
\end{list}
\end{definition}
\begin{remarque}
The morphism $A\rightarrow B$ is $hf$ (resp $hfp$) if and only if the morphism $q_{c}A\rightarrow q_{c}B$ is $hf$ (resp $hfp$). The morphism $q_{c}B\rightarrow B$ is always $hf$.  
\end{remarque}

\begin{lemme}\label{stabfc}
The hf (resp hfp) morphisms are stable under composition.
\end{lemme}    
Proof\\The proofs for $hf$ morphisms and $hfp$ morphisms are analogous so let us prove it for $hf$ morphisms. Let $A\rightarrow B\rightarrow C$ be the composition of two hf morphisms. There is a diagram
\begin{center}
$\xymatrix{q_{c}A\ar[r]\ar[d]& q_{c}B\ar[r]\ar[d]&q_{c}C\ar[d]\\A\ar[r]&B\ar[r]&C}$
\end{center}
and forgetful functors $F_{1}:sq_{c}C-mod\rightarrow sq_{c}B-mod$ and $F_{2}:sq_{c}B-mod\rightarrow sq_{c}A-mod$. The image $F_{1}(q_{c}C)$ of $q_{c}C$ is homotopically finitely presented in $sq_{c}B-mod$ hence weakly equivalent to a retract of a finite homotopical colimit of $q_{c}B$ in $Ho(sq_{c}B-mod)$. The forgetful functor $F_{2}$ preserves retracts, equivalences, finite colimits, cofibrant objects and cofibrations whose domain is cofibrant. Thus it also preserves finite homotopical colimit and sends $q_{c}C$ to a retract of a finite homotocipal colimit of $q_{c}B$ in $Ho(sq_{c}A-mod)$. As $q_{c}B$ is homotopically finitely presented in $sq_{c}A-mod$, and as homotopically finitely presented objects are stable under retracts, equivalences and finite homotopical colimit, $C$ is sent by $F_{2}\circ F_{1}$ in $sq_{c}A-mod_{c}$. Hence $A\rightarrow C$ is finite.
\vskip 2pt
$\hfill\blacklozenge$

\begin{lemme}\label{stabfhp}
The hf (resp hfp) morphims are stable under homotopic pushout of simplicial monoids. 
\end{lemme}
Proof:\\The proofs for $hf$ morphisms and $hfp$ morphisms are analogous so let us prove it for $hf$ morphisms. Let $A\rightarrow B$ and $A\rightarrow C$ be in $sComm(\mathcal{C})$ such that the first is finite. Let $q_{cA}$ be the cofibrant replacement of $q_{c}A-alg$, it is weakly equivalent to $q_{c}$ and the object $q_{c_{A}}B$ is homotopically finitely presented in $sq_{c}A-mod$. Let us prove that $B\otimes_{A}^{h}C\backsimeq q_{cA}B\otimes_{q_{c}A} q_{c}C$ (in $Ho(q_{c}A-mod)$, Reedy lemma) is homotopically finitely presented in $q_{c}C-mod$. The forgetful functor $sq_{c}C-mod\rightarrow sq_{c}A-mod$ preserves filtered colimits and weak equivalences hence it preserves homotopical filtered colimits. Thus the derived functor $-\otimes_{q_{c}A} q_{c}C$ preserves homotopically finitely presented objects. So $B\otimes_{A}^{h}C$ is homotopically finitely presented in $Ho(q_{c}C)$.
\vskip 2pt
$\hfill\blacklozenge$

\subsection{A Definition for Smoothness}
\begin{definition}\label{fs}
A morphism $A\rightarrow B$ in $sComm(\mathcal{C})$ is formally smooth if the morphism $B\otimes_{A}^{h}B\rightarrow B$ is $hf$. 
\end{definition}
\begin{remarque}
This definition does not generalise the definition of formal smoothness in the sense of rings. However, the corresponding notion of smoothness is a generalisation of the classical notion of smoothness, as it will be proved in this article.
\end{remarque}

\begin{proposition}\label{stabc}
Formally smooth morphisms are stable under composition.
\end{proposition}
Proof:\\by previous remarks, it can be assumed that $A$ is cofibrant in $sComm(\mathcal{C})$, $B$ is cofibrant in $sA-alg$ and $C$ is cofibrant in $sB-alg$. Let $A\rightarrow B\rightarrow C$ be the composition of two formally smooth morphisms . The morphisms $B\coprod_{A}B\rightarrow B$ and $C\coprod_{B}C\rightarrow C$ are $hf$. The following diagram commutes and is clearly cocartesian :
\begin{center}
$\xymatrix{C\\C\coprod_{A}C\ar[u]\ar[r]&C\coprod_{B}C\ar[lu]\\B\coprod_{A}B\ar[u]\ar[r]&B\ar[u]}$
\end{center}
Thus, if it is cofibrant for the Reedy stucture, it will be homotopically cocartesian. The morphisms $B\backsimeq B\coprod_{A} A\rightarrow B\coprod_{A}B$ and $B\coprod_{A}B\rightarrow C\coprod_{A}C$ are images by the left Quillen functor $colim$, of clear Reedy cofibrations (see \cite{a} for a descriptions of these cofibrations), thus are cofibrations. In particular $B\coprod_{A}B$ and $C\coprod_{A}C$ are cofibrant and the diagram considered is Reedy cofibrant. Finally, the morphism $C\coprod_{A}C\rightarrow C\coprod_{B}C$ is $hf$ as a pushout of $hf$ morphisms and  $C\coprod_{A}C\rightarrow C$ is $hf$ as a composition of $hf$ morphisms.
\vskip 2pt
$\hfill\blacklozenge$

\begin{proposition}\label{stabhp}
Formally smooth morphism are stable under homotopic pushout.
\end{proposition}
\noindent{Proof:\\Let $u: A\rightarrow B$ be a formally smooth morphism and $C$ be a commutative $A$-algebra. By previous remarks it can be assumed that $A$ and $c$ are cofibrants in $sComm(\mathcal{C})$ and that $B$ is cofibrant in $sA-alg$. Let $D$ denote the homotopic pushout of $B\otimes_{A}C$ and $u'$ denote the morphism from $B$ to $D$. Clearly:} 
\begin{center}
$D\otimes_{C}D\backsimeq B\otimes_{A}C\otimes_{C}B\otimes_{A}C\backsimeq B\otimes_{A}D$
\end{center}
Thus the following diagram commutes :
\begin{center}
$\xymatrix{&B\otimes_{A}B\ar[ld]\ar[d]^{Id\otimes_{A}f}\ar[r]^{m_{B}}&B\ar[d]^{f}\\D\otimes_{C}D\ar[r]&B\otimes_{A}D\ar[r]_{m_{D}}&D}$
\end{center}
And is cocartesian : 
\begin{center}
$B\otimes_{B\otimes_{A}B}B\otimes_{A}B\otimes_{A}C\backsimeq B\otimes_{A}C \backsimeq D$
\end{center}
Moreover it is clearly cofibrant as $B\otimes_{A}-$ preserve cofibrations. Finally by stability of $hf$ morphism under homotopic pushouts, the morphism $C\rightarrow D$ is formally smooth.
\vskip 2pt
$\hfill\blacklozenge$

\begin{definition}
Let $A$ be in $sComm(\mathcal{C})$ and $M$ be in $sA-mod$.
\begin{enumerate}
\item The object $M$ is $n$-truncated if $Map_{s\mathcal{C}}(X,M)$ is $n$-truncated in $sSet$, $\forall\;X\in s\mathcal{C}$.
\item The Tor-Dimension of $M$ in $sA-mod$ is defined by
\begin{center}
$Tordim_{A}(M)=inf\{n\;st\;M\otimes^{h}_{A}X\;is\;n+p-truncated\;\forall\;X\in\;sA-mod\;p-truncated\}$
\end{center}
\item A morphism of monoids $A\rightarrow B$ has \textit{Tor dimension $n$} if $Tordim_{A}(B)=n$.
\end{enumerate}
\end{definition}
\begin{lemme}\label{stabtd}
\textit{Tor dimension zero} morphisms are stable under composition and homotopic pushout.
\end{lemme}
Proof:\\Let $A\rightarrow B\rightarrow C$ be the composition of two \textit{Tor dimension zero} morphisms. Let $M$ be a $p$ truncated $A$-module,
\begin{center}
$M\otimes_{A}C\backsimeq M\otimes_{A}B\otimes_{B}C$.
\end{center}
As $Tordim_{A}(B)=0$, $M\otimes_{A}B$ is a $p$ truncated $B$-module. As $Tordim_{B}(C)=0$, $M\otimes_{A}C$ is a $p$ truncated $C$-module. Thus $Tordim_{A}(C)=0$.\vskip 2pt
Let $A\rightarrow B$ be a \textit{Tor dimension zero} morphism and $A\rightarrow C$ be a morphism in $Comm(\mathcal{C})$. Let $M$ be in $C$-mod and let $D$ denote the pushout $B\otimes_{A}C$. We have
\begin{center}
$M\otimes_{C}B\otimes_{A}C\backsimeq M\otimes_{A}B$.
\end{center} 
Thus, $Tordim_{A}(B)=0$ implies $Tordim_{C}(D)=0$. 
\vskip 2pt
$\hfill\blacklozenge$
\begin{definition}
A morphism $A\rightarrow B$ in $Comm(\mathcal{C})$ is smooth if it is formally smooth, $hfp$ and has Tor-Dimension zero. A morphism of affine scheme is smooth if the corresponding morphism of monoids is smooth. We say that an affine scheme $X$ is smooth if the morphism $X\rightarrow Spec(1)$ is smooth. 
\end{definition}
\begin{theoreme}
Smooth morphisms are stable under composition and homotopic pushout.
\end{theoreme}
Proof:\\
This a a corollary of \ref{stabtd}, \ref{stabc}, \ref{stabhp}, \ref{stabfc} and \ref{stabfhp}.
\vskip 2pt
$\hfill\blacklozenge$

\section{Simplicial Presheaves Cohomology}
In the article \cite{t1}, B. Toën define a cohomology for a connected and pointed simplicial presheaf. We will define here a cohomology for a general simplicial presheaf. This theory will be used to fin examples of smooth morphisms of commutative monoids (in sets). The references cited in this section are \cite{t1}, \cite{gj} and \cite{j}.
\subsection{Definitions}
In this section, $\mathcal{D}$ is a category and $sPr(\mathcal{D})$ is the category of simplicial presheaves over $\mathcal{D}$.
\begin{definition}(\cite{gj}$VI.3$)\\
Soit $F\in sPr(\mathcal{D})$. The tower of $n$-truncations of $F$ is a Postnikov tower:
\begin{center}
$\xymatrix{...\ar[r]&\tau_{\leq n}F\ar[r]&\tau_{\leq n-1}F\ar[r]&...\ar[r]&\tau_{\leq 1}F\ar[r]&\tau_{\leq 0}F}$.
\end{center}   
\end{definition} 
\begin{definition}
Let $F$ be a simplicial presheaf. 
\begin{list}{$\triangleright$}{}
\item The functor $\pi_{0}(F):\mathcal{D}\rightarrow Ens$ is defined by $\pi_{0}(F):X\rightarrow \pi_{0}(F(X))$.
\item The category $(\mathcal{D}/F)_{0}$ is the full subcategory of $sPr(\mathcal{D})/F$ whose objects are in $\mathcal{D}$.
\item The functor $\pi_{n}(F):(\mathcal{D}/F)_{0}\rightarrow Ens$ is defined by $\pi_{n}(F)(X,u)=\pi_{n}(F(X),u)$.
\end{list}
\end{definition}
\begin{definition}
Let $G$ be a simplicial group.
\begin{list}{$\triangleright$}{}
\item The bisimplicial set $E(G,1)$ is defined by $E(G,1)_{p,q}=G_{p}^{q}$.
\item The classifying space of $G$, denoted $K(G,1)$, is given by the diagonal of the bisimplicial set $E(G,1)/G$. More precisely $K(G,1)_{n}=G_{n}^{n}/G_{n}$. It is abelian if $G$ is abelian.
\item The endofunctor of abelien groups $K(G,1)^{\circ n}$ is denoted $K(G,n)$.
\end{list}
\end{definition}
\begin{remarques}
As the diagonal of $E(G,1)$ is pointed (by identity), the simplicial set $K(G,1)$ is also pointed. In particular, $\pi_{n}(K(G,1),*)\backsimeq \pi_{n-1}(G,e_{g})$. This construction is functorial (in $G$) and then extends to presheaves of simplicial groups.
\end{remarques}
\subsection{Simplicial presheaves Cohomology}
It is necessary to work in the proper category to construct a cohomology for a simplicial presheaf $F$ which is not connected or pointed. In fact the $1$-truncation of $F$ is the nerve $NG$ of a groupoid $G$ and in the category $sPr(\mathcal{D})/NG$, $F$ becomes connected and pointed. But in this category, there is no clear construction fro classifying spaces. The solution of this problem is given by a Quillen equivalence with the category $sPr(\mathcal{D}/G)$, for a well chosen category $\mathcal{D}/G$. We choose now a simplicial presheaf $F$.
\subsubsection{The Category of Presheaves}
\paragraph{The left adjoint functor $\widetilde{(-)}$}
\begin{definition}
The category $\mathcal{D}/G$ is the category whose objects are couples $(X,x)$, $x:X\rightarrow NG$, and whose morphisms from $(X,x)$ to $(Y,y)$ are couples $(f,u)$ where $f:X\rightarrow Y$ and $u:y\circ f\backsimeq x$ in $G(X)\backsimeq \pi_{1}F(X)$. 
\end{definition}
Next step is to construct a functor  $\widetilde{(-)}:\mathcal{D}/G\rightarrow sPr(\mathcal{D}/NG)$.
\begin{definition}
Let $(X,x,)$ be in $\mathcal{D}/G$. Define a presheaf of groupoïds $G_{X,x}$ on $\mathcal{D}$. Ths image of $S\in \mathcal{D}$ is the groupoid described as follow
\begin{list}{-}{}
\item The objects are triples $(u,y,h)$, $u:S\rightarrow X$, $y\in G(S)$, and $h:x\circ u\backsimeq y\in G(S)$. 
\item A morphism from $(u,y,h)$ to $(u',y',h')$ is an endomorphism $k$ of $S$ such that $k^{*}(h:x\circ u\rightarrow y)=h':x'\circ u'\rightarrow y'$.
\end{list}
Let $\breve{X}$ denote the nerve of this groupoid.
\end{definition}
\begin{remarque}
There is a commutative diagram of presheaves of groupoids
\begin{center}
$\xymatrix{X\ar[rr]^{x}\ar[rd]_{j}&&G\\&G_{X,x}\ar[ru]_{l}&}$
\end{center}
where $l$ is the projection on $G$ and $j$ is given for $S\in\mathcal{D}$ by $j(S):u\in Hom_{\mathcal{D}}(S,X)\rightarrow (u,x\circ u,Id)\in G_{X,x}$.\vskip 2pt
\noindent Applying the functor nerve, one get a morphism $\breve{x}:=Nl:\breve{X}\rightarrow G$. It defines a functor  
\begin{center}
$\breve{(-)}:\mathcal{D}/G\rightarrow sPr(\mathcal{D})/NG$\vskip 2pt
$(X,x)\rightarrow (\breve{X},\breve{x})$
\end{center}
\end{remarque}
\begin{definition}
The functor $\widetilde{(-)}:\mathcal{D}/G\rightarrow sPr(\mathcal{D})/NG$ is defined by
\begin{center}
$\widetilde{(-)}:(X,x)\rightarrow (\widetilde{X},\widetilde{x}):=Q(\breve{X},\breve{x})$
\end{center}
where $Q$ is a cofibrant replacement in $sPr(\mathcal{D})/NG$.
\end{definition}
\begin{remarques}
This functor has a kan extension to $sPr(\mathcal{D}/G)$, still denoted
\begin{center}
$\widetilde{(-)}:sPr(\mathcal{D}/G)\rightarrow Spr(\mathcal{D})/NG$.
\end{center}
In facts, the category $sPr(\mathcal{D}/G)$ is equivalent to the category $sPr(\mathcal{D})^{NG}$ defined in \cite{j} and the equivalence of category we are constructing is constructed in a different way and a more general situation in \cite{j}.  
\end{remarques}
\paragraph{The right adjoint functor $(-)_{1}$}\vskip 1pt
\noindent We construct now the (right) adjoint of $\widetilde{(-)}$, denoted $(-)_{1}$.
\begin{definition}
The functor $(-)_{1}:sPr(\mathcal{D})/NG\rightarrow sPr(\mathcal{D}/G)$ is defined by
\begin{center}
$(-)_{1}:(H,u)\rightarrow H_{1}:=(X,x)\rightarrow \underline{Hom}^{\Delta}_{sPr(\mathcal{D})/NG}((\widetilde{X},\widetilde{x}),(H,u))$
\end{center}
\end{definition}
where $\underline{Hom}^{\Delta}$ is the simplicial $Hom$. As $(\widetilde{X},\widetilde{x})$ is constructed cofibrant, the functor $(-)_{1}$ is right Quillen and its adjoint is then left Quillen. We prove now that $R(-)_{1}$ commute with homotopy colimits. We need to recall first some properties.
\begin{definition}
Let $(H,h)$ be in $sPr(\mathcal{D})/NG$ and $(X,x)$ be in $\mathcal{D}/G$. Define an object $(H_{X},h_{x})$ by the homotopy pullback diagram
\begin{center}
$\xymatrix{H_{X}\ar[r]\ar[d]&H\ar[d]^{h}\\X\ar[r]_{x}&NG}$
\end{center}
\end{definition}
\begin{lemme}
Let $(H,f)$ be an homotopy colimit, $H\backsimeq Hocolim(H_{i})$, in $sPr(\mathcal{D})/NG$ and let $(X,x)$ be in $\mathcal{D}/G$.
\begin{list}{$\diamond$}{}
\item There is an isomorphism $H_{X}\backsimeq Hocolim (H_{i})_{X}$ in $Ho(sPr(\mathcal{D})/NG)$.
\item There is an isomorphism $RH_{1}(X)\backsimeq Map_{sPr(\mathcal{D})/X}((X,Id),(H_{X},h_{x}))$.
\end{list}
\end{lemme}
\begin{corollaire}
The functor $R(-)_{1}$ commute with homotopy colimits.
\end{corollaire}
Proof:\\
Let $H$ be isomorphic to $Hocolim(H_{i})$ and $(X,x)$ be in $\mathcal{D}/G$.
\begin{center}
$RH_{1}(X)\backsimeq Map_{sPr(\mathcal{D})/X}((X,Id),([Hocolim(H_{i})]_{X},[Hocolim(h_{i})]_{x}))$\vskip 3pt
$\backsimeq Map_{sPr(\mathcal{D})/X}((X,Id),(Hocolim[(H_{i})_{X}],Hocolim[(h_{i})_{x}]))\backsimeq Hocolim(R(H_{i})_{1}(X))$
\end{center}
\vskip 2pt
$\hfill\blacklozenge$
\paragraph{The Equivalence}
\begin{proposition}
The Quillen functors $\widetilde{(-)}$ and $(-)_{1}$ define a Quillen equivalence.
\end{proposition}
Proof:\\
The functor $\widetilde{(-)}$ commutes with homotopy colimits and as any object in $sPr(\mathcal{D}/G)$ is an homotopy colimit of representable objects $H\backsimeq Hocolim(X_{i})$, its image can be computed in terms of representable objects, i.e. $\widetilde{H}\backsimeq Hocolim(\widetilde{X_{i}})$. The short exact sequence $H_{1}\rightarrow H\rightarrow \tau_{\leq 1}H$ proves that $(-)_{1}$ preserves weak equivalences. Then 
\begin{center}
$(\widetilde{H})_{1}\backsimeq (Hocolim(\widetilde{X_{i}}))_{1}\backsimeq Hocolim(X_{i})\backsimeq H$.
\end{center}
If $H$ is cofibrant in $sPr(\mathcal{D}/G)$ and $H'$ is fibrant in $sPr(\mathcal{D})/NG$, we consider a morphism between short exact sequences 
\begin{center}
$\xymatrix{H\ar[r]\ar[d]&\widetilde{H}\ar[d]\ar[r]&NG\ar[d]^{Id}\\H'_{1}\ar[r]&H'\ar[r]&NG}$
\end{center}
Applying the functors $\pi_{i}$, it is clear that if $\widetilde{H}\rightarrow H'$ is an equivalence, so is $H\rightarrow H'_{1}$ and reciprocally, if $H\rightarrow H'_{1}$, the homotopic fibers of $\widetilde{H}\rightarrow H'$ upon $NG$ are equivalences thus so is $\widetilde{H}\rightarrow H'$.
\vskip 2pt
$\hfill\blacklozenge$
\subsubsection{The Cohomology}
\begin{definition}
Let $F$ be in $sPr(\mathcal{D})$, a local system on $F$ is a presheaf of abelian groups on $\mathcal{D}/G$, where $G$ verifies $NG\backsimeq \tau_{\leq1}F$. A Morphism of local system is a morphism of presheaves of abelian groups. The category of local systems on $F$ will be denoted $sysloc(F)$. The n-th classifying space of $M$ is denoted $K(M,n)$ and its image by $L\widetilde{(-)}$ is denoted $L\tilde{K}(M,n)$.
\end{definition}
\begin{remarque}
The object $L\tilde{K}(M,n)$ is characterised up to equivalence by the fact that $\pi_{n}(L\tilde{K}(M,n))\backsimeq M$, $\pi_{1}(L\tilde{K}(M,n))\backsimeq \pi_{1}(F)$, $\pi_{0}(L\tilde{K}(M,n))\backsimeq \pi_{0}(F)$ and that its other homotopy presheaves of groups are trivial.
\end{remarque}
\begin{definition}
Let $F$ be in $sPr(\mathcal{D})$ and $M$ be a local system on $F$. The n-th cohomology group of $F$ with coefficient in $M$ is 
\begin{center}
$H^{n}(F,M):=\pi_{0}Map_{sPr(\mathcal{D})/NG}(F,L\tilde{K}(M,n))$
\end{center}
\end{definition}
The standard example of local system is $\pi_{n}$. Indeed, it has been defined on $(\mathcal{D}/F)_{0}$ but it clearly lifts to $\mathcal{D}/G$.\vskip 2pt
The important theorem is here.
\begin{theoreme}\label{thmimp}
Let $G$ be a grouoid. For all $m$, the functor
\begin{center}
$H^{m}(NG,-):Sysloc(NG)\rightarrow Ab$\vskip 2pt
$M\rightarrow H^{m}(NG,M)$
\end{center}
is isomorphic to the n-th derived functor of the functor $H^{0}(NG,-)$.
\end{theoreme}
Proof:\\
There is an equivalence between the category of simplicial abelian group presheaves, denoted $sAb(\mathcal{D}/G)$, on $\mathcal{D}/G$ and the category of complexes of abelian group presheaves with negative or zero degree, denoted $C^{-}(\mathcal{D}/G,Ab)$. This is a generalisation of Dold-Kan correspondence. There is a correspondence between quasi-isomorphism of complexes and weak equivalences of simplicial presheaves, and then an induced equivalence between the homotopical categories :
\begin{center}
$\Gamma: D^{-}(\mathcal{D}/G,Ab)\backsimeq Ho(sAb(\mathcal{D})/G)$
\end{center}
The derived functors of $H^{0}$ are then given by 
\begin{center}
$H^{m}_{der}(\mathcal{D}/G,M)\backsimeq Hom_{D^{-}(\mathcal{D}/G,Ab)}(\mathbb{Z},M[m])$
\end{center}
Where $\mathbb{Z}$ is regarded as a complex concentrated in degree zero and $M[m]$ is concentrated in degree $-m$, with value $M$. As $\Gamma(\mathbb{Z})$ is the constant presheaf with fiber $\mathbb{Z}$, still denoted $\mathbb{Z}$, and as $\Gamma(M[m])$ is equivalent to $K(M,m)$, $\Gamma$ induces an isomorphism :
\begin{center}
$Hom_{D^{-}(\mathcal{D}/G,Ab)}(\mathbb{Z},M[m])\backsimeq Hom_{Ho(sAb(\mathcal{D}/G))}(\mathbb{Z},K(M,m))$
\end{center}
Finally, the adjunction between the abelianisation functor,denoted $\mathbb{Z}(-)$ from $sPr(\mathcal{D}/G)$ to $sAb(\mathcal{D}/G)$ and the forgetful functor gives
\begin{center}
$Hom_{D^{-}(\mathcal{D}/G,Ab)}(\mathbb{Z},M[m])\backsimeq Hom_{Ho(sPr(\mathcal{D}/G))}(*,K(M,m))\backsimeq H^{m}(NG,M)$.
\end{center}  
\vskip 2pt
$\hfill\blacklozenge$
\subsubsection{Obstruction Theory}
There is an homotopic pullback diagram in $sPr(\mathcal{D}/G)$:
\begin{center}
$\xymatrix{\tau_{\leq n}F_{1}\ar[r]\ar[d]&\txt{*}\ar[d]\\\tau_{\leq n-1}F_{1}\ar[r]&K(\pi_{n}(F),n+1)}$
\end{center}
As $F_{1}$ is $1$-connex, this pullback diagram is a (functorial) generalisation to presheaf of the diagram given by the proposition 5.1 of \cite{gj}. By the quillen equivalence $((-)_{1},\widetilde{(-)})$, there is an homotopic pullback diagram:
\begin{center}
$\xymatrix{\tau_{\leq n}F\ar[r]\ar[d]&NG\ar[d]\\\tau_{\leq n-1}F\ar[r]&L\widetilde{K}(\pi_{n}(F),n+1)}$
\end{center}
If $H\rightarrow \tau_{\leq n-1}F$ is a morphism in $Ho(sPr(\mathcal{D})/NG)$, it has a lift to $\tau_{\leq n}F$ if and only if it is send to a zero element in the group
\begin{center}
$\pi_{0}Map_{sPr(\mathcal{D})/NG}(H,L\widetilde{K}(\pi_{n}(F),n+1))$
\end{center}
This group can be described in terms of cohomology. Indeed, if $G'$ is a groupoid such that $NG'\backsimeq \tau_{\leq 1}H$. Let $u$ denote the morphism $u:NG'\rightarrow NG$. To simplify the notations, we still write $H$ for what we sould call $u^{*}H$. There is a Quillen adjunction:
\begin{center}
$\xymatrix{sPr(\mathcal{D})/NG'\ar@<1ex>[r]^{u^{*}}&\mathcal{D}/NG\ar@<1ex>[l]^{-\times_{NG}NG'}}$
\end{center}
which induces an isomorphism
\begin{center}
$Map_{sPr(\mathcal{D})/NG}(H,L\widetilde{K}(\pi_{n}(F),n+1))\backsimeq Map_{sPr(\mathcal{D})/NG'}(H,L\widetilde{K}(\pi_{n}(F),n+1)\times_{NG}NG')$.
\end{center}
There is a clear weak equivalence $L\widetilde{K}(\pi_{n}(F),n+1)\times_{NG}NG'\backsimeq L\widetilde{K}(\pi_{n}(F)\circ u^{*},n+1)$, thus :
\begin{center}
$\pi_{0}Map_{sPr(\mathcal{D})/NG}(H,L\widetilde{K}(\pi_{n}(F),n+1))\backsimeq H^{n+1}(H,\pi_{n}(F)\circ u^{*})$
\end{center}
\subsection{Simplicial Modules Cohomology}
It is well known that for a commutative monoid $B$ in $(Set,\times, \mathbb{F}_{1})$, there is an equivalence 
\begin{center}
$sPr(\mathcal{B}B)\backsimeq sB-mod$
\end{center}
where $\mathcal{B}B$ is the category with one object with a set of endomorphisms isomorphic to $B$. We will identify these two categories in this part. Let now $A$ be a commutative monoid in sets and $B\rightarrow A$ be a morphism of commutative monoids. We are in a particular case of previous section, the category $\mathcal{D}$ is $\mathcal{B}B$ and the presheaf of groupoids $G$ is just $A$. Let $M$ be a local system on $\mathcal{B}B$, there is an isomorphism
\begin{center}
$H^{n}(A,M)\backsimeq \pi_{0}Map_{sB-mod/A}(A,L\widetilde{K}(M,n+1))$.
\end{center}
Let $Z$ denote the abelianization functor from $B-mod/A$ to the category of abelian group objects in $B-mod/A$, denoted $Ab(B-mod/A)$. There is an equivalence between $Ab(B-mod/A)$ and the category of $A$ graduated $Z(B)$-modules, denoted $Z(B)-mod^{A-grad}$. The following functor realizes this equivalence, its inverse is the forgetful functor.
\begin{center}
$\Theta:(\xymatrix{M\ar[r]^{f}&A})\in Ab(sB-mod/A)\rightarrow \oplus_{m\in A}f^{-1}(m)\in Z(B)-mod^{A-grad}$\vskip 5pt
\end{center}
This equivalence lifts to simplicial categories and it is easy to see that
\begin{eqnarray}\label{for}
H^{n+1}(A,M)\backsimeq \pi_{0}Map_{Z(B)-mod^{A-grad}}(Z(A),L\widetilde{K}(M,n+1))
\end{eqnarray}
Here is the proposition that interrests us.
\begin{proposition}\label{bl}
Let $B\rightarrow A$ be a morphism of commutative monoids in sets. The morphism $B\rightarrow A$ is $hf$ if and only if
\begin{list}{$\diamond$}{}
\item $Z(A)$ is homotopically finitely presented in $Z(B)-mod^{A-grad}$.
\item $A$ is homotopically finitely presented for the $1$-truncated model structure i.e. in the category $B-Gpd$.
\end{list}
\end{proposition}
Proof\\
Let us prove first the easiest part. Let $A$ be an homotopically finitely presented object in $sB-mod$. Let $sB-mod^{\leq 1}$ denotes the category $sB-mod$ endowed with its $1$-truncated model structure. In the adjunctions 
\begin{center}
$\xymatrix{sB-mod\ar@<1ex>[r]^{Id}&sB-mod^{\leq 1}\ar@<1ex>[l]^{Id}}$
\vskip 5pt
$\xymatrix{sB-mod\ar@<1ex>[r]{Z}&sZ(B)-mod/A\ar@<1ex>[l]^{i}}$
\vskip 5pt
$\xymatrix{sB-mod/A\ar@<1ex>[r]^{Z}{Z}&sZ(B)-mod^{A-grad}\ar@<1ex>[l]^{i}}$
\end{center}
the left adjoint functors preserve weak equivalences and cofibrations thus the right adjoints preserve homotopically finitely presentable objects.\vskip 4pt
Let us now prove the hardest part. We start with this lemma:
\begin{lemme}\label{mz}
There exists $m_{0}\in \mathbb{N}$ such that for any local system $M$ and all $n\geq m_{0}$
\begin{center}
$H^{n}(A,M)\backsimeq *$.
\end{center}
\end{lemme}
Proof\\
The isomorphism \ref{for} proves that the cohomology of $A$ is isomorphic to the Ext functors of $Z(A)$ in $sZ(B)-mod^{A-grad}$. Moreover, there is an equivalence of abelian categories
\begin{center}
$sZ(B)-mod^{A-grad}\backsimeq C^{-}(\mathcal{B}B/A,Ab)$
\end{center}
which induces by \ref{thmimp} an equivalence with the derived functors of $H^{0}$. In particular as $Z(A)$ is homotopically finitely presented, the derived functors of $H^{0}$ vanished after a set rank denoted $m_{0}$.
\vskip 2pt
$\hfill\blacklozenge$
\begin{remarque}\label{ses}
Two corollaries comes now. They are a consequence of this lemma and the following short exact sequence, $C\in A/sB-mod$
\begin{center}
$\xymatrix{Map_{sB-mod/\tau_{\leq n-1}C}(A,\tau_{\leq n}C)\ar[r]& Map_{sB-mod}(A,\tau_{\leq n}C) \ar[d]\\Map_{sB-mod/NG}(A,L\widetilde{K}(\pi_{n}(C),n+1))& Map_{sB-mod}(A,\tau_{\leq n-1}C)\ar[l]}$
\end{center}
\end{remarque}
\begin{corollaire}\label{iso}
Let $\xymatrix{A\ar[r]^{v}&C}$ be in  $A/sB-mod$. For all $i\geq 1$, for all $n\geq n_{i}=n_{0}+i+1$ 
\begin{center}
$\pi_{0}Map_{sB-mod}(A,\tau_{\leq n-1}C)\backsimeq \pi_{0}Map_{sB-mod}(A,\tau_{\leq n}C)$
\vskip 4pt
$\pi_{i}(Map_{sB-mod}(A,\tau_{\leq n-1}C),v)\backsimeq \pi_{i}(Map_{sB-mod}(A,\tau_{\leq n}C),v)$
\end{center} 
\end{corollaire}
Proof\\
We first prove that the simplicial set $Map_{sB-mod/\tau_{\leq n-1}C}(A,\tau_{\leq n}C)$ is not empty. There are pushout squares :
\begin{center}
$\xymatrix{A\times^{h}_{L\widetilde{K}(\pi_{n}(C),n+1)}NG\ar[r]\ar[d]&\tau_{\leq n}C\ar[r]\ar[d]&NG\ar[d]_{s}\\A\ar[r]&\tau_{\tau_{\leq n-1}C}\ar[r]&L\widetilde{K}(\pi_{n}(C),n+1)\ar@/_1pc/[u]_{p}}$
\end{center}
where $p\circ s= Id$. There are then equivalences
\begin{center}
$Map_{sB-mod/\tau_{\leq n-1}F}(A,\tau_{\leq n}C)\backsimeq Map_{sB-mod/L\widetilde{K}(\pi_{n}(C),n+1)}(A,NG)\backsimeq Map_{sB-mod/A}(A,A\times^{h}_{L\widetilde{K}(\pi_{n}(C),n+1)}NG)$.
\end{center}
Let $f$ be the morphism from $A$ to $L\widetilde{K}(\pi_{n}(C),n+1)$. There is a morphism $p\circ f:A\rightarrow NG$. As the cohomology of $A$ vanished for $n\geq n_{0}$, the elements $s\circ p\circ f$ and $f$ of the cohomology group are equals and thus 
\begin{center}
$p\circ f\in \pi_{0}Map_{sB-mod/L\widetilde{K}(\pi_{n}(C),n+1)}(A,NG)$.
\end{center}
Then, for $i=0$, the corollary is a clear consequence of lemma \ref{mz} and the short exact sequence of remark \ref{ses}.\\ 
Now, Let us study the case $i>0$. As $NG\times^{h}_{L\widetilde{K}(\pi_{n}(C),n+1)}NG\backsimeq L\widetilde{K}(\pi_{n}(C),n+1)$, we obtain 
\begin{center}
$A\times^{h}_{L\widetilde{K}(\pi_{n}(C),n+1)}NG\backsimeq L\widetilde{K}(\pi_{n}(C)\circ v^{*},n)$
\end{center}
Thus
\begin{center}
$\pi_{i}(Map_{sB-mod/\tau_{\leq n-1}C}(A,\tau_{\leq n}C),v)\backsimeq \pi_{i}((Map_{sB-mod/A}(A,L\widetilde{K}(\pi_{n}(C)\circ v^{*},n)),q)\backsimeq H^{n-i}(A,\pi_{n}(C))$
\end{center}
where $q$ is the natural morphism from $A$ to $L\widetilde{K}(\pi_{n}(C)\circ v^{*},n)$. We deduce then the result from lemma \ref{mz} and the short exact sequence of remark \ref{ses}.  
\vskip 2pt
$\hfill\blacklozenge$
\begin{corollaire}
Let $\xymatrix{A\ar[r]^{v}&C}$ be in $A/sB-mod$. The pointed tower of fibrations
\begin{center}
$(Map_{sB-mod}(A,\tau_{\leq n}C),v)$
\end{center} 
converges completly in the sense of \cite{gj}.
\end{corollaire}
Proof\\
It can be checked with the corollary 2.21 of the complete convergence lemma of \cite{gj}.
\vskip 2pt
$\hfill\blacklozenge$
\begin{corollaire}\label{tqs}
For all $i\geq 0$, all $n\geq n_{i}$ and all $\xymatrix{A\ar[r]^{v}&C}$ in $A/sB-mod$, there are isomorphisms
\begin{center}
$\pi_{i}(Map_{sB-mod}(A,C),v)\backsimeq lim_{n\in \mathbb{N}}\pi_{i}(Map_{sB-mod}(A,\tau_{\leq n}C),v)\backsimeq \pi_{i}(Map_{sB-mod}(A,\tau_{\leq n_{i}}C),v)$
\end{center}
\end{corollaire}
Proof\\
The first isomorphism is a consequence of Milnor exact sequence (\cite{gj}, 2.15) and the vanishing of the $lim^{1}$ induced by the complete convergence. The second isomorphism is a consequence of corollary \ref{iso}
\vskip 2pt
$\hfill\blacklozenge$
\vskip 2pt
Let us now recall a well known lemma with which we will prove the last technical lemma necessary for the proof of \ref{bl}.
\begin{lemme}\label{tsvf}
Let 
\begin{center}
$\xymatrix{X\ar[r]^{f}\ar[d]&Y\ar[d]\\Z\ar[r]_{g}&T}$
\end{center}
be a commutative square in sSet where $g$ is a weak equivalence. The morphism $f$ is a weak equivalence if and only if for all $z\in Z$, the homotopic fibers $X_{z}$ and $Y_{g(z)}$ are simultaneously empty and equivalent when not empty.
\end{lemme}
Here is the last technical lemma:
\begin{lemme}\label{hopfini}
Let $C\backsimeq Hocolim_{\alpha\in \Theta}(C_{\alpha})$ be an homotopical filtered colimit. There is a weak equivalence in $sSet$
\begin{center}
$Map_{sB-mod}(A,C)\backsimeq Hocolim Map_{sB-mod}(A,C_{\alpha})$.
\end{center}
\end{lemme}
Proof\\
By induction on the truncation level $n$ of $C$. This is an hypothesis of \ref{bl} for $n=1$. Let us assume that is is true for $n-1$. Let $C$ be an $n$-truncated object in $sB-mod$ and $\bar{u}$ be in $Hocolim Map_{sB-mod}(A,\tau_{n-1}C_{\alpha})$, represented by $u\in Map_{sB-mod}(A,\tau_{n-1}C_{\alpha_{0}})$. Let $\tilde{u}$ denote its image in $Map_{sB-mod}(A,C)$. The filtered hocolimit along $\Theta$ is weak equivalent to the hocolimit along $\alpha_{0}/\theta$. We will use previous lemma, computing the fibers along $u$ as in the following diagram:
\begin{center}
$\xymatrix{Hocolim_{\alpha_{0}/\Theta}Map_{sB-mod/\tau_{n-1}C_{\alpha}}((A,u_{\alpha}),C_{\alpha})\ar[r]\ar[d]&Map_{sB-mod/\tau_{n-1}C}((A,\tilde{u}_{n-1}),C)\ar[d]\\Hocolim_{\alpha_{0}/\Theta}Map_{sB-mod}(A,C_{\alpha})\ar[r]\ar[d]&Map_{sB-mod}(A,C)\ar[d]\\Hocolim_{\alpha_{0}/\Theta}Map_{sB-mod}(A,\tau_{n-1}C_{\alpha})\ar[r]&Map_{sB-mod}(A,\tau_{n-1}C)}$
\end{center}   
Where $u_{\alpha}:\xymatrix{A\ar[r]^{u}&C_{\alpha_{0}}\ar[r]&C_{\alpha}}$ and $\tilde{u}_{n-1}:\xymatrix{A\ar[r]^{\tilde{u}}&C\ar[r]&\tau_{n-1}C}$.\vskip 4pt
Let us show first that the fibers are simultaneously empty. The naturel morphism 
\begin{center}
$Hocolim_{\alpha_{0}/\Theta}Map_{sB-mod}(A,\tau_{n-1}C_{\alpha})\rightarrow Map_{sB-mod}(A,\tau_{n-1}C)$\vskip 4pt
$\bar{u}\rightarrow \tilde{u}$
\end{center}
induces the naturel morphism on cohomology groups
\begin{center}
$Hocolim_{\alpha_{0}/\Theta}H^{n+1}(A,\pi_{n}C_{\alpha})\rightarrow H^{n+1}(A,\pi_{n}C)$
\end{center}
which is a weak equivalence. Indeed, the $H^{n}$ are isomorphic to $Ext$ functors in $sZ(B)-mod^{A-grad}$ which commute with filtered hocolimits by the first hypothesis of \ref{bl}. The images of $\bar{u}$ and $\tilde{u}$ in the cohomology groups vanish then simultaneously, and the fibers are simultaneously empty. \vskip 4pt
Let us assume now that the fibers are unempty and prove that they are equivalent. The functors $\pi_{i}$ commute with homotopical filtered colimits, applying them on the fibers, we get the following natural morphism
\begin{center}
$colim_{\alpha_{0}/\Theta}\pi_{i}Map_{sB-mod/\tau_{n-1}C_{\alpha}}((A,u_{\alpha}),C_{\alpha})\rightarrow \pi_{i}Map_{sB-mod/\tau_{n-1}C}((A,\tilde{u}_{n-1}),C)$
\end{center} 
As these $\pi_{i}$ are in fact isomorphic to $H^{n-1}$, these morphisms are isomorphisms. By \ref{tsvf}, this ends the proof of the lemma.
\vskip 2pt
$\hfill\blacklozenge$
\vskip 6pt
Let us now prove \ref{bl}.\\
Let $v:A\rightarrow C$ be in $A/sB-mod$ such that $C\backsimeq Hoclolim(C_{\alpha})$. Let us prove that the morphism
\begin{center}
$Hocolim(Map_{sB-mod}(A,C_{\alpha}))\rightarrow Map_{sB-mod}(A,C)$
\end{center}
is a weak equivalence. Let $i$ be a positive integer, to check if the image of this morphism by $\pi_{i}$ is an isomorphism, we can just consider the case $C$ $n$-truncated by \ref{tqs}. As the truncation commuta with homotopical filtered colimits, this is a consequence of \ref{hopfini}. This ends the proof of \ref{bl}.
\vskip 2pt  
$\hfill\blacklozenge$
  
\section{Examples}
\subsection{The Category $(\mathbb{Z}-mod,\otimes_{\mathbb{Z}},\mathbb{Z})$}
In classical algebraic geometry, the notion of (projective) resolution is obtained using chain complex of modules or rings. In facts, considering the correspondence of Dold-Kan this method is equivalent to taking cofibrant resolution in the simplical category (cf \cite{q}).
\begin{theoreme}
(Dold-Kahn correspondance)\\Let $A$ be a ring. There is an equivalence of categories:
\begin{center}
$sA-mod\backsimeq Ch(A-mod)^{\geq0}$  and $\;\forall i\;\pi_{i}(Map(\mathbb{Z},X))\backsimeq H_{i}(X)$.
\end{center}
In particular, it induces a correspondence between weak equivalences and quasi-isomorphisms.
\end{theoreme}

\begin{remarque}Let $A$ be a ring. Generating cofibrations of $Ch(A-mod)^{\geq0}$ are levelwise equal to $\{0\}\rightarrow A$ or $Id_{A}$.
\end{remarque}
\begin{definition}Let $A$ be a rings,  $M,N$ be two $A$-modules. 
\begin{enumerate}
\item Define $Tor_{*}^{A}(M,N):= H_{*}(M\otimes^{L}_{A}N)$. 
\item Define $Ext_{A}^{*}(M,N):=H^{*}(R\underline{Hom}_{A-mod}(M,N))$.
\item Define the projective dimension of $M$ by:  
\begin{center}
$ProjDim_{A}(M):=inf\{n\;st\;Ext_{A}^{n+1}(M,-)=\{0\}\}$. 
\end{center}
\item Define the Tor-dimension of $M$ by: 
\begin{center}
$TorDim_{A}(M):= inf\{n\;st\;\forall\;X\;p-truncated\;Tor^{A}_{i}(M,X)=\{0\}\;\forall i>n+p\}$
\end{center}
\end{enumerate}
\end{definition}
\begin{remarque}
The functor of Dold-Kan correspondence is a strong monoidal functor, as a consequence the $Tor$ dimension can be computed with $\pi_{i}$ instead of $H_{i}$.
\end{remarque}

\begin{lemme}Let $X$ be in $Ho(sSet)$ and $M$ be in $s\mathbb{Z}-mod$ (resp $sA-mod$, for $A$ a ring)
\begin{list}{$\diamond$}{}
\item The object $X$ is $n$-truncated if and only if $Map(*,X)\backsimeq Map(S^{i},X)$ $\forall\;i>n$ in $Ho(sSet)$. 
\item The object $M$ is $n$-truncated if and only if $Map_{s\mathbb{Z}-mod}(\mathbb{Z},M)$ (resp $Map_{sA-mod}(A,Z)$) is $n$-truncated in $Ho(sSet)$.
\end{list}
\end{lemme}
Proof\\
For the first statement, by\ref{tsvf}, we can consider equivalently the homotopic fibers of this morphism upon $Map(*,X)$. The fiber of $Map(*,X)$ is a point and the fiber of $Map(S^{i},X)$ is $Map_{sSet/*}(S^{i},X)$. As $\pi_{j}Map_{sSet/*}(S^{i},X)\backsimeq \pi_{i+j}(X)$, the equivalence is clear.\vskip 4pt
For the second statement, any object in $s\mathbb{Z}-mod$ is an homotopical colimit of free objects, i.e. $\forall$ $N\in s\mathbb{Z}-mod$ there exists a family of sets $(\lambda_{i})_{i\in I}$ such that $qN\backsimeq hocolim_{I}\coprod_{\lambda_{i}}\mathbb{Z}$ in $Ho(s\mathbb{Z}-mod)$ . Assume that $Map_{s\mathbb{Z}-mod}(\mathbb{Z},M)$ is $n$-truncated. $Map_{s\mathbb{Z}-mod}(N,M)\backsimeq holim_{I}\prod_{\lambda_{i}}(Map_{s\mathbb{Z}-mod}(\mathbb{Z},M))$, hence is an homotopical limit of $n$-truncated objects. by $i$, $n$-truncated objects in $sSet$ are clearly stable under homotopical limits.
\vskip 2pt
$\hfill\blacklozenge$

\begin{lemme}(cf \cite{tv})
Let $u:A\rightarrow B$ be in $s\mathbb{Z}-mod$. The morphism $u$ is flat if and only if 
\begin{enumerate}
\item The natural morphism $\pi_{*}(A)\otimes_{\pi_{0}(A)}\pi_{0}(B)\rightarrow \pi_{0}(B)$ is an isomorphism. 
\item The morphism $\pi_{0}(u)$ is flat.
\end{enumerate}
In particular, if $A$ is cofibrant and $n$-truncated, $u$ flat implies $B$ $n$-truncated.
\end{lemme}
\begin{remarque}\label{flat}\cite{tv}
Let $A\rightarrow B$ be in $\mathbb{Z}-alg$. The morphism $A\rightarrow B$ is flat if and only if $TorDim_{A}(B)=0$.
\end{remarque}

We give now the lemmas necessary to the theorem of comparison of the notions of smoothness in rings and relative smoothness. 

\begin{lemme}\label{noet}
Let $A\rightarrow B$ be a smooth morphism of rings. There exists a pushout square
\begin{center}
$\xymatrix{A'\ar[r]\ar[d]&B'\ar[d]\\A\ar[r]&B}$
\end{center}
such that $A'\rightarrow B'$ is a smooth morphism of noetherian rings. 
\end{lemme}
Proof:\\This is the affine case in the corollary $17.7.9(b)$ of \cite{ega4}.
\vskip 2pt
$\hfill\blacklozenge$

\begin{lemme}
Let $A\rightarrow B$ and $A\rightarrow C$ be two morphisms in $\mathbb{Z}-alg$. If $B$ is a perfect complex of $B\otimes_{A}B$ modules then $D:=B\otimes_{A}C$ is a perfect complex of $D\otimes_{C}D$ modules.  
\end{lemme}
\noindent Proof:\\Perfect complexes are clearly stable under base change. As $D\otimes_{C}D\backsimeq B\otimes_{A}D$, the natural morphism $D\otimes_{C}D\rightarrow D$ is a pushout of $B\otimes_{A}B\rightarrow B$ hence $D$ is a perfect complex.
\vskip 2pt
$\hfill\blacklozenge$

\begin{lemme}
Let $A$ be a noetherian ring. Every flat $A$-module of finite type is projective.
\end{lemme}

\begin{lemme}\label{dim}
Assume that $A$ is a noetherian ring and consider $A\rightarrow B\in \mathbb{Z}-alg$, $B$ of finite type. There is an equivalence between 
\begin{enumerate}
\item The ring $B$ is of finite Tor-dimension on $A$.
\item The ring $B$ is of finite projective dimension on $A$.
\end{enumerate}
\end{lemme}
The part $ii\Rightarrow i$ is clear, if $B$ has a finite projective resolution $0\rightarrow P_{n}\rightarrow... \rightarrow B$, then for $i\geq n$, $Tor^{i+1}(M,-)\backsimeq Tor^{i-n}(P_{n+1},-)$ and $P_{n+1}=0$.\vskip 3pt Reciprocally, if $TorDim_{A}b<+\infty$, let $...\rightarrow P_{n}\rightarrow...\rightarrow B$ be a free resolution of $B$. The module $P_{n}/im(P_{n+1})$ has Tor dimension $0$ by previous formula hence is flat by \ref{flat}. As $A$ is noetherian and $B$ is of finite type, it is projective and we have a clear finite projective resolution.\vskip 2pt
$\hfill\blacklozenge$

\begin{lemme}\label{afi}
Let $u:A\rightarrow B$ be in rings. Assume that $A$ is an algebraically closed field, then there is an equivalence 
\begin{list}{$\diamond$}
\item The morphism $u$ is formally smooth in the sense of rings.
\item Any morphism $x:B\rightarrow A$ in rings provides $A$ with a structure of $B$-module of finite projective dimension over $B$.
\end{list}
\end{lemme}
\begin{lemme}\label{alf}
Let $u:A\rightarrow B$ be a finitely presented flat morphism in rings. The morphism $u$ is smooth if and only if for all algebraically closed field $K$ under $A$, $K\rightarrow K\otimes_{A}B$ is smooth.
\end{lemme}

\begin{theoreme}
A morphism $A\rightarrow B$ in $\mathbb{Z}-alg$ is smooth in the sense of rings if and only if
\begin{enumerate}
\item The ring  $B$ is finitely presented in $A-alg$.
\item The morphism $A\rightarrow B$ is flat.
\item The ring $B$ is a perfect complex of $B\otimes_{A}B$-modules.
\end{enumerate}
\end{theoreme}
Proof:\\
Let us now prove the first part of the theorem. Assume that $A\rightarrow B$ is smooth. $i$ and $ii$ are clear.\vskip 1pt
Let us prove $iii$. By \ref{noet}, as $iii$ is stable under pushout, we just have to prove it for $A$ and $B$ noetherian. Let us prove first that $B\otimes_{A}B\rightarrow B$ is of finite $Tor$ dimension (hence of finite projective dimension by \ref{dim}).\\Let $L$ be an algebraically closed field in $A-alg$. Set $B_{L}:= B\otimes_{A}L$. Clearly
\begin{center}
$B\otimes_{B\otimes_{A}B}L\backsimeq B_{L}\otimes_{B_{L}\otimes_{L} B_{L}}L$
\end{center} 
hence computing the Tor dimension of $B$ over $B\otimes_{A}B$ is equivalent to compute the Tor dimension of $B_{L}$ over $B_{L}\otimes_{L}B_{L}$. The morphism $L\rightarrow B_{L}\rightarrow B_{L}\otimes_{L}B_{L}$ is smooth, by composition of smooth morphisms, over an algebraically closed field. The ring $B_{L}\otimes_{L}B_{L}$ is then smooth on afield, hence regular. Now, $B_{L}$ is a module of finite type on this regular ring thus it is a perfect complex on it. In particular, it is of finite projective dimension hence of finite $Tor$ dimension. Finally, $B$ is of finite $Tor$ dimension hence of finite projective dimension over $B\otimes_{A}B$. As previously, $B$ of finite type over $B\otimes_{A}B$. As these rings are noetherian, $B$ is a perfect complex. Indeed, $B$ has a finite projective resolution by $(P_{i})$. Each $P_{i}$ is of finite $Tor$ dimension hence of finite projective dimension.\vskip 6pt

Let us prove the second part of the theorem. Let $A\rightarrow B$ be a morphism of rings verifying $i$, $ii$ and $iii$. Let $K$ be an algebraically closed field under $A$. We will use \ref{alf} and \ref{afi}.\\
Let $x:B\rightarrow K$ be in $\mathbb{Z}-alg$. The following commutative diagram is an homotopic pushout:
\begin{center}
$\xymatrix{B\otimes_{K}B\ar[r]^{Id\otimes_{K}x}\ar[d]&B\otimes_{K}K\backsimeq B\ar[d]^{x}\\B\ar[r]_{x}&K}$
\end{center}
Thus $K$ has finite projective dimension in $B-mod$. Finally, by \ref{alf}, $K\rightarrow B$ is smooth in the sense of rings. As it is true for any $K$, by \ref{afi}, $A\rightarrow B$ is smooth in the sense of rings.  
\vskip 2pt
$\hfill\blacklozenge$
\vskip 2pt
Here is now the comparison theorem.
\begin{theoreme}
Let $A\rightarrow B$ be a morphism of rings. It is smooth if and only if it is smooth in the sense of rings. 
\end{theoreme}
Proof\\ The two following lemmas, and remark \ref{flat} prove the theorem.
\begin{lemme}\cite{tv}
Let $A\rightarrow B$ be a morphism in $\mathbb{Z}-alg$. 
\begin{enumerate}
\item if $A\rightarrow B$ is $hfp$, then it is finitely presented in $\mathbb{Z}-alg$.
\item if $A\rightarrow B$ is smooth and finitely presented, then it is $hfp$.
\end{enumerate}
\end{lemme}
\begin{lemme}ref{tv}
Let $A\rightarrow B$ be a morphism of rings. The ring $B$ is a perfect complex of $B$-modules if and only if $A\rightarrow B$ is $hf$.
\end{lemme}
\vskip 2pt
$\hfill\blacklozenge$

\subsection{The category $Set$}

The most difficult problem consists in finding examples of formally smooth morphisms. The Lemma \ref{bl} gives us a characterisation of these morphisms in the relative context $\mathcal{C}=Set$.\\
The functor nerve and the functor "fundamental groupoid" define a Quillen equivalence between the category $sB-mod$ endowed with its $1$-truncated model structure and the category $B-Gpd$. Moreover, this last category is compactly generated and thus its filtered $Hocolim$ can be computer as filtered colimits. Here is the formula to do this
\begin{lemme}
Let $\mathcal{I}$ be a filtered diagram and $F:\mathcal{I}\rightarrow Gpd$. The colimit of $F$ consists of
\begin{list}{$\diamond$}{}
\item On objects 
\begin{center}
$(Colim F)_{0}:=Colim(fg\circ F)$
\end{center}
where $fg$ is the forgetful functor from $Gpd$ to $Set$.
\item On morphisms, for $\bar{x},\bar{y}\in Colim(fg\circ F)$ represented by $x\in F(i)$ and $y\in F(i')$.  There exists $k$ under $i$ and $i'$ such that
\begin{center}
$Hom_{Hocolim(F)}(\bar{x},\bar{y}):=Colim_{k/\mathcal{I}}(Hom_{F(j)}((l_{i,j})_{*})(x),(l_{i',j})_{*})(y))$
\end{center}
where $l_{i,j}:i\rightarrow j$ and $l_{i',j}:i'\rightarrow j$.
\end{list}
\end{lemme} 
We also need to describe the derived enriched Homs. 
\begin{lemme}
Let $B$ be a monoid in $Set$. There is an equivalence of categories between $Ho(B-Gpd)$ and the category [$B-Gpd]$ whose objects are $B$-groupoids and morphisms are isomorphism classes of functors. In particular, for two $B$-groupoids $G$ and $G'$,  $R\underline{Hom}_{B-gpd}^{\Delta\leq1}(G,G')\backsimeq \underline{Hom}_{[B-gpd]}^{\Delta\leq1}(G,G')$ in $Ho(Gpd)$, where the exponent $\Delta\leq1$ means that the Homs are enriched on groupoids.
\end{lemme}

\begin{lemme}
The commutative monoid $\mathbb{N}$ is homotopically finitely presented for the $1$-truncated model structure i.e. in the category $(\mathbb{N}\times \mathbb{N})-Gpd$.
\end{lemme}
Let $\mathbb{N}²$ denotes $\mathbb{N}\times \mathbb{N}$.  Let $F:\mathcal{J}\rightarrow Gpd$ be a funnctor from a filtered diagram $\mathcal{I}$ to $Gpd$. We have to prove
\begin{center}
$Hocolim(\underline{Hom}^{\Delta\leq1}_{[\mathbb{N}²-Gpd]}(\mathbb{N},F(-)))\backsimeq \underline{Hom}^{\Delta\leq1}_{[\mathbb{N}-Gpd]}(\mathbb{N},Hocolim(F))$
\end{center}
We let the reader verify that the following functor denoted $\varphi$ define an equivalence of groupoids.\vskip 2pt
Let $\bar{H}$ be in $Hocolim(\underline{Hom}^{\Delta\leq1}_{[\mathbb{N}²-Gpd]}(\mathbb{N},F(-)))$ represented by $H\in \underline{Hom}_{[\mathbb{N}²-gpd]}(\mathbb{N},F(j))$. We define $\varphi$ on objects by
\begin{center}
$\varphi:\bar{H}\rightarrow \hat{H}:=n\rightarrow \bar{H(n)}$
\end{center}
Now, by construction,  any morphism $\bar{\eta}$ in $Hocolim(\underline{Hom}^{\Delta\leq1}_{[\mathbb{N}²-Gpd]}(\mathbb{N},F(-)))$ has a representant $\eta:G\rightarrow G'\in Hom_{\underline{Hom}_{[\mathbb{N}²-gpd]}(\mathbb{N},F(j))}(G,G')$. We define $\varphi$ on morphisms by
\begin{center}
$\varphi:\bar{\eta}\rightarrow \hat{\eta}:=n\rightarrow \bar{\eta_{n}}$
\end{center}
\vskip 2pt
$\hfill\blacklozenge$
\begin{lemme}
The commutative group $\mathbb{Z}$ is homotopically finitely presented for the $1$-truncated model structure i.e. in the category $(\mathbb{Z}\times \mathbb{Z})-Gpd$.
\end{lemme}
Proof\\This is the same proof as previous lemma, replacing $\mathbb{N}$ by $\mathbb{Z}$.
\vskip 2pt
$\hfill\blacklozenge$
\begin{corollaire}
The morphisms $\mathbb{F}_{1}\rightarrow \mathbb{N}$ and $\mathbb{F}_{1}\rightarrow \mathbb{Z}$ are smooth. In particular, the affine scheme $Gl_{1,\mathbb{F}_{1}}\backsimeq Spec(\mathbb{Z})$, also denoted $\mathbb{G}_{m,\mathbb{F}_{1}}$ in \cite{tva}, is smooth.
\end{corollaire} 
Proof\\
They are clearly $hfp$ and of $Tor$ dimension zero. Their diagonal is hf for the $1$-truncated model structure, thus, we just have to check that the diagonal of their abelianisation is $hf$ in the simplicial graduated category given in \ref{bl}. The abelianisation of $\mathbb{N}$ is $\mathbb{Z}[X]$ and the abelianisation of $\mathbb{Z}$ is $\mathbb{Z}(X)$, and the morphisms $\mathbb{Z}[X]\otimes_{\mathbb{Z}}\mathbb{Z}[X]\rightarrow \mathbb{Z}[X]$ and $\mathbb{Z}(X)\otimes_{\mathbb{Z}}\mathbb{Z}(X)\rightarrow \mathbb{Z}(X)$ are $hf$ respectively in $s(\mathbb{Z}[X]\otimes_{\mathbb{Z}}\mathbb{Z}[X])-Mod^{\mathbb{N}-grad}$ and $s(\mathbb{Z}(X)\otimes_{\mathbb{Z}}\mathbb{Z}(X))-Mod^{\mathbb{Z}-grad}$.
\vskip 2pt
$\hfill\blacklozenge$
\begin{corollaire}
Pour tout $n$, le schéma $Gl_{n,\mathbb{F}_{1}}$ est lisse.
\end{corollaire}
Proof\\This scheme is isomorphic to $Spec(\prod_{E_{n}}\coprod_{E_{n}}\mathbb{Z})$ (\cite{tva}), where $E_{n}$ is the set of integers from $1$ to $n$, thus as coproducts in $Comm(Set)$ are products in $Set$, it is isomorphic to $Spec(\mathbb{Z}^{n²})$ . The product in set is the tensor product, thus as a finite tensor product of finite colimits of homotopically finitely presentable object, this monoid is homotopically finitely presentable, i.e. $\mathbb{F}_{1}\rightarrow \mathbb{Z}^{n²}$ is a morphism $hfp$. For the same reason,  the Tor dimension is still zero. We need then to prove that a finite tensorisation of the formally smooth morphism $*\rightarrow \mathbb{Z}$ (in the relative sense \ref{fs}) by itself is still formally smooth. The pushout diagram
\begin{center}
$\xymatrix{\mathbb{Z}^{2}\ar[r]\ar[d]&\mathbb{Z}\ar[d]\\\mathbb{Z}^{k}\ar[r]&\mathbb{Z}^{k-1}}$
\end{center}
proves that $\mathbb{Z}^{k}\rightarrow \mathbb{Z}^{k-1}$ is $hf$ for any integer $k$ and by composition $\mathbb{Z}^{2k}\rightarrow\mathbb{Z}^{k}$ is $hf$ for any integer $k$. Finaly for every $n$, $\mathbb{F}_{1}\rightarrow \mathbb{Z}^{n²}$ is smooth, hence $Gl_{n},\mathbb{F}_{1}$ is smooth.
\vskip 2pt
$\hfill\blacklozenge$

\subsection{Some Others examples}
If ($\mathcal{C},\otimes,1$) is a symmetric monoidal category as described in the preliminaries, its associated category of simplcial objects has simplicial Homs,denoted $\underline{Hom}^{\Delta}$, and there is an adjunction
\begin{center}
$\xymatrix{s\mathcal{C}\ar@<1ex>[r]^{\underline{Hom}^{\Delta}(1,-)}&sSet\ar@<1ex>[l]^{sK_{0}}}$
\end{center}
where $sK_{0}((X_{n})_{n\in \mathbb{N}})=(\coprod_{X_{n}}1)_{n\in \mathbb{N}}$. One verifies easily that as $1$ is cofibrant, finitely presentable, and as $\underline{Hom}^{\Delta}(1,-)$ preserves weak equivalences (by construction of the model structure on $\mathcal{C}$), the functor $sK_{0}$ preserve homotopically finitely presentable objects. In particular, $sK_{0}$ preserves $hf$ morphisms and formally smooth morphisms. Restricting the adjunction to the categories of algebra, where weak equivalences and homotopical filtered colimits are obtained with the forgetful functor, it is also clear that $sK_{0}(u)$ preserves $hfp$ morphisms. We write then the following proposition.
\begin{proposition}
Let $u:A\rightarrow B$ be a smooth morphism in $Comm(Set)$, then $sK_{0}(u)$ is smooth if and only if $sK_{0}(B)$ is of finite $Tor$ dimension over $sK_{0}(A)$.
\end{proposition}
This gives particular examples. Indeed, in every context the affine line correspond to the morphism $1\rightarrow 1[X]:=\coprod_{\mathbb{N}}1$ and the scheme $\mathbb{G}_{m}$ to the morphism $1\rightarrow 1(X):=\coprod_{\mathbb{Z}}1$. We write then the following theorem.
\begin{theoreme}
The affine line and the scheme $\mathbb{G}_{m}$ are smooth in any context where, respectively, $1[X]$ and $1(X)$ are of finite $Tor$ dimension over $1$. 
\end{theoreme}
This theorem can be applied in particular to the context $\mathbb{N}-mod$. The following lemma provides us, in this context, examples of morphisms of Tor-dimension $0$.
\begin{lemme}
Let $A\rightarrow B$ be in $Comm(\mathbb{N}-mod)$ such that $B$ is free over $A$. The monoid $B$ has Tor-dimension $0$ over $A$.
\end{lemme}
Proof :\\Let $M\in A-mod$ be a $n$-truncated module. There exists a set $\lambda$ such that $B\backsimeq \coprod_{\lambda}A$. Thus $B\otimes_{A}^{L}M'\backsimeq Coprod_{\lambda}QM$ in $Q_{c}A-mod$ where $Q,Q_{c}$ are cofibrant replacement respectively in $Q_{c}A-mod$ and $Comm(\mathbb{N}-mod)$. Thus as this coproduct is a product in set, we get
\begin{center}
$B\otimes_{A}^{L}M'\backsimeq Colim_{\lambda'fini \subset \lambda}\prod_{\lambda'}QM$
\end{center}
As functors $\pi_{i}$ commute with products in sets and filtered colimits, the $Tor$ dimenson of $B$ over $A$ is zero.
\begin{theoreme}Examples in $\mathbb{N}-mod$.
\begin{list}{$\diamond$}{}
\item The affine line in $\mathbb{N}-mod$, $\mathbb{A}^{1}_{\mathbb{N}}$, is smooth. 
\item The scheme $\mathbb{G}_{m,\mathbb{N}}$ relative to $\mathbb{N}-mod$ is smooth.
\end{list}
\end{theoreme}
We conclude with a last theorem
\begin{theoreme}
Let $\mathcal{C}$ be a relative context in the sense of \cite{m} and $A\rightarrow B$ be a Zariski open immersion in $Comm(\mathcal{C})$, with $A$ cofibrant in $Comm(\mathcal{C})$ and $B$ cofibrant in $A-alg$. The morphism $A\rightarrow B$ is smooth.
\end{theoreme}
Proof\\
A Zariski open immersion is always formally smooth, its diagonal is even an isomorphism. Thus we will need to prove that it is $hfp$ and of $Tor$ dimension zero. First, if there exists $f\in A_{0}$,an object of the underlying set of $A$, such that $B\backsimeq A_{f}$, the result is clear. Indeed, $A_{f}$ is given by a filtered colimit of $A$ thus is of $Tor$ dimension zero. Let us prove that it is $hfp$.  It is clear that $A\rightarrow A[X]$ is homotopically finitely presented, then as everything is cofibrant, we can write $A_{f}$ as a finite colimit of $A[X]$ (\cite{m}) which is in facts a finite homotopical colimit and thus finally  $A\rightarrow A_{f}$ is $hfp$.\vskip 2pt
Now if $B$ define a Zariski open object of $A$, we can write $B$ it as a cokernel of products of $A_{f}$. As functors $\pi_{i}$ commute with products, the products preserve weak equivalences and it is then clear that  $A\rightarrow B$ is $hfp$. For the $Tor$ dimension, recall that there is a finite family of functor reflecting isomorphisms $B-mod\rightarrow A_{f}-mod$. Let $M$ be a $p$-truncated $A$-module. This family sends $M\otimes^{L}_{A}B$ and its $n$-truncations,$n>p$ to the same module $QM_{f}$ (Q is the cofibrant replacement of $A-mod$) thus clearly $M\otimes^{L}_{A}B$ is $p$ truncated and $TorDim_{A}(B)=0$.
\vskip 2pt
$\hfill\blacklozenge$

\end{document}